\documentclass[11pt]{amsart}

\input xy
\xyoption{all}

\begin{document}

\date{September 15, 2003}

\title[Homology of Invariant Group Chains]{The Homology of Invariant Group Chains}

\author{Kevin P.~Knudson}\thanks{Partially supported by NSF grant no.~DMS-0242906.}
\address{Department of Mathematics and Statistics, Mississippi State University,
P.~O.~Drawer MA, Mississippi State, MS 39762}
\email{knudson@math.msstate.edu}

\newtheorem{theorem}{Theorem}[section]
\newtheorem{prop}[theorem]{Proposition}
\newtheorem{lemma}[theorem]{Lemma}
\newtheorem{cor}[theorem]{Corollary}
\newtheorem{conj}[theorem]{Conjecture}
\newtheorem{definition}[theorem]{Definition}
\newtheorem{remark}[theorem]{Remark}

\newcommand{\zz}{{\mathbb Z}}
\newcommand{\zq}{{\mathbb Q}}
\newcommand{\gm}{{\mathbb G}_m}
\newcommand{\ra}{\rightarrow}
\newcommand{\lra}{\longrightarrow}
\newcommand{\bop}{\bigoplus}
\newcommand{\zl}{{\mathbb Z}/\ell}
\newcommand{\zn}{{\mathbb Z}/n}
\newcommand{\rbar}{\overline{R}}
\newcommand{\kbar}{\overline{k}}
\newcommand{\zr}{{\mathbb R}}
\newcommand{\zc}{{\mathbb C}}
\newcommand{\abar}{\underbar{A}}
\newcommand{\cA}{{\mathcal A}}
\newcommand{\cQ}{{\mathcal Q}}
\newcommand{\cbar}{\underbar{C}}
\newcommand{\cH}{{\mathcal H}}

\maketitle

\section*{Introduction}
Let $Q$ be a finite group acting on a group $G$ as a group of automorphisms. This action induces
an action of $Q$ on the standard bar complex $C_\bullet(G)$ for computing the homology
of $G$; denote the subcomplex of invariants by $C_\bullet(G)^Q$.  In this paper we study
the homology of the complex $C_\bullet(G)^Q$.  The resulting homology groups are denoted
$H_\bullet^Q(G)$.

Before summarizing the results, a brief history is in order.  A quick literature search
and consultation with experts in equivariant topology revealed nothing.  Indeed, everyone
I asked was more or less alarmed that I wanted to know about invariants of a {\em chain}
complex (instead of a cochain complex).  So I forged ahead with the cumbersome calculations
using the bar complex and obtained a few results.  For example, if $A$ is an abelian group
with trivial $G$- and $Q$-actions, then we can define homology groups with coefficients
$H_\bullet^Q(G;A)$.  The inclusion of complexes $C_\bullet(G;A)^Q\ra C_\bullet(G;A)$ induces
a homomorphism $H_\bullet^Q(G;A)\ra H_\bullet(G;A)$ whose image clearly lies in the subgroup
$H_\bullet(G;A)^Q$ of $Q$-invariant homology classes.  The following is easily proved.

\medskip

\noindent {\bf Proposition \ref{invariantiso}.}  {\em Suppose that $|Q|$ is invertible in $A$.  Then the natural map
$$H_\bullet^Q(G;A)\lra H_\bullet(G;A)^Q$$ is an isomorphism.}

\medskip

In certain cases, we can construct transfer maps
$$H_\bullet^Q(G;A)\lra H_\bullet^Q(K;A)$$
for $Q$-stable subgroups $K$, satisfying the usual properties.  One then obtains
the following expected result.

\medskip

\noindent {\bf Proposition \ref{orderkill}.} {\em If $G$ is finite, then $H_i^Q(G;A)$ is annihilated
by $|G|$ for all $i>0$.}

\medskip

The existence of transfers also allows the calculation of $H_\bullet^{\zz/2}(\zn)$ for certain
$n$, where $\zz/2$ acts on $\zn$ by $x\mapsto -x$.

\medskip

\noindent {\bf Theorem \ref{nodd}.}  {\em If $n\not\equiv 0 \mod 4$, then
$$H_\bullet^{\zz/2}(\zn;\zz) \cong H_\bullet(\zn;\zz)^{\zz/2}.$$}

The case $n\equiv 0\mod 4$ presented difficulty, however.  None of the techniques used to prove
Theorem \ref{nodd} work in this case.  Moreover, it is easy to show that $H_1^{\zz/2}(\zn) \cong
\zz/2\oplus\zz/2$ in this case, and computer calculations of the low-dimensional groups for $n=4$
indicated that something fairly complicated was going on.  So I put this project on the back
burner in hopes that a new idea would present itself.

The primary motivation for studying the $H_\bullet^Q(G)$ arose from joint work with Mark Walker
\cite{knudsonwalker}.  There, we were confronted with calculating the groups $H_\bullet(BG(\rbar)/\Gamma)$,
where $R$ is a nice affine domain over a field $k$, $\rbar$ its absolute integral closure, $G$ an
algebraic group over $k$, and $\Gamma$ the Galois group $\text{Gal}(\overline{k(R)}/k(R))$.  The
calculation is practically impossible in general, but we were able to handle $R=\zr$, $\rbar=\zc$,
$G=\gm$, and $\Gamma = \text{Gal}(\zc/\zr) = \zz/2$.  For some time, I worked under the mistaken
assumption that the complex $C_\bullet(G(\rbar))^\Gamma$ computed this homology (there is an isomorphism
of the groups of chains in each degree, but it does not commute with the boundary map).  However,
the work in \cite{knudsonwalker} did provide some insight into the problem considered here.

Consider the norm map 
$N:C_\bullet(BG/Q;A)\lra C_\bullet(G;A)^Q$
defined by
$$\sigma\otimes a \stackrel{N}{\mapsto} \sum_{q\in Q} q\sigma\otimes a.$$
This provides a map
$$H_\bullet(BG/Q;A)\lra H_\bullet^Q(G;A)$$ whose composition with the obvious
map 
$$H_\bullet^Q(G;A)\lra H_\bullet(G;A)^Q$$
is the transfer map encountered in the study of the homology of quotient spaces.
The homology of the quotient complex $$D_\bullet=C_\bullet(G;A)^Q/C_\bullet(BG/Q;A)$$
is computable in some cases.  We have the following results.

\medskip

\noindent {\bf Theorem \ref{hiz}.} {\em Let $\zz/2$ act on $\zz$ via $n\mapsto -n$.  Then for $i>0$, we have
$$H_i^{\zz/2}(\zz;\zz) \cong \begin{cases}
                               \zz/2 & i\;\text{\em odd} \\
                               0  & i\;\text{\em even}.
                               \end{cases}$$}

\medskip

One might have thought that $H_i^{\zz/2}(\zz)$ would vanish for $i\ge 2$ since $B\zz\simeq S^1$.
The problem arises, however, because $S^1$ is not a $\zz/2$-equivariant model for $B\zz$.

If $n\equiv 0 \mod 4$, the calculation of $H_\bullet^{\zz/2}(\zn)$ is possible using techniques in
equivariant topology.  One first reduces to $n=2^s$, $s\ge 2$.  

\medskip

\noindent {\bf Theorem \ref{n0mod4}.} {\em Assume $s\ge 2$.  Then for all $\ell>0$,
$$H_\ell^{\zz/2}(\zz/2^s;\zz) \cong \begin{cases}
                                           (\zz/2)^{k} & \ell=2k \\
                                           (\zz/2)^{2k} & \ell=4k-3 \\
                                           \zz/2^s\oplus (\zz/2)^{2k} & \ell=4k-1.
                                        \end{cases}$$}

This paper is organized as follows.  In Section \ref{definition} we develop the basic properties
of $H_\bullet^Q(G;A)$ and the cohomological version $H^\bullet_Q(G;A)$.  Section \ref{bredon} relates
our construction to the homology of quotient spaces.  In Section \ref{basic} we make some basic calculations.  Section
\ref{transfer} is devoted to the transfer maps and their application.  Finally, in Section \ref{calculate}
we perform the above-mentioned calculations.  The approach to calculations is to be as low-tech as
possible for as long as possible; that is, we try to work with chains directly, resorting to more
sophisticated techniques only when absolutely necessary.

\medskip

\noindent {\em Acknowledgements.}  I would like to thank Paul Goerss, John Greenlees, and Mark Walker for
valuable discussions about this material.  I also thank those who have listened to my talks about this
for their suggestions; these include Eric Friedlander, Rick Jardine, and 
Stewart Priddy.

\section{Basic Constructions}\label{definition}  Denote by $C_\bullet(G)$ the bar complex for computing the
homology of the group $G$:
$$C_i(G) = \zz\{[g_1|\cdots |g_i]: g_j\in G\}$$
\begin{eqnarray*}
\lefteqn{d([g_1|\cdots |g_i]) =}  & & \\
   &  &  [g_2|\cdots |g_i] + \sum_{k=1}^{i-1} (-1)^k[g_1|\cdots |g_kg_{k+1}|\cdots |g_i] 
+ (-1)^i[g_1|\cdots |g_{i-1}]
\end{eqnarray*}
(that is, $C_\bullet(G) = B_\bullet\otimes_G \zz$, where $B_\bullet$ is the bar resolution).  If $A$ is
an abelian group, set $C_\bullet(G;A) = C_\bullet(G)\otimes A$ with boundary $d\otimes 1$.  Let $Q$ be a
finite group acting as a group of automorphisms of $G$.  The action of $Q$ on $G$ induces an action
on $C_\bullet(G;A)$ by
$$q([g_1|\cdots |g_i]\otimes a) = [q(g_1)|\cdots |q(g_i)]\otimes a.$$
Denote the subcomplex of invariant chains by $C_\bullet(G;A)^Q$.

\begin{definition} $$H_\bullet^Q(G;A) = h_\bullet(C_\bullet(G;A)^Q)$$
$$H^\bullet_Q(G;A) = h^\bullet(\text{\em Hom}(C_\bullet(G)^Q,A)).$$
\end{definition}

First note the following basic fact.

\begin{lemma}\label{maps}  Denote by $G^Q$ the subgroup of $Q$-invariant elements of $G$.  Then there are
natural maps
$$H_\bullet(G^Q;A)\lra H_\bullet^Q(G;A)\lra H_\bullet(G;A)^Q.$$
\end{lemma}

\begin{proof}  Note that there are inclusions of complexes
$$C_\bullet(G^Q;A)\lra C_\bullet(G;A)^Q\lra C_\bullet(G;A)$$
and these induce maps
$$H_\bullet(G^Q;A)\stackrel{f_*}{\lra} H_\bullet^Q(G;A)\stackrel{i_*}{\lra} H_\bullet(G;A).$$
Note that the image of $i_*$ lies in the subgroup $H_\bullet(G;A)^Q$ of $Q$-invariants since a
$Q$-invariant cycle $z\in C_\bullet(G;A)^Q$ gives rise to an invariant homology class
$[z]\in H_\bullet(G;A)^Q$.
\end{proof}

\section{Homology of quotient spaces}\label{bredon}
Consider the CW-complex $BG$ whose chain complex is $C_\bullet(G)$.  The $Q$-action on 
$G$ induces a cellular action on $BG$ with the property that if an element of $Q$ fixes
a cell, then it fixes it pointwise.  This gives $BG$ the structure of a $Q$-CW-complex
and the quotient space $BG/Q$ is also a CW-complex.

Observe that we have
\begin{eqnarray*}
C_n(BG/Q) & = & \zz\{[g_1|\cdots |g_n]: g_i\in G\}/\langle [q(g_1)|\cdots |q(g_n)] = [g_1|\cdots |g_n]\rangle \\
    & = & C_n(BG)_Q
\end{eqnarray*} 
with the boundary map being induced by the boundary map in $C_\bullet(G)$.  Recall that for any $Q$-module
$M$, there is a norm map $M_Q\ra M^Q$ whose kernel and cokernel are annihilated by $|Q|$.  This map is
natural for maps $M\ra M'$ of $Q$-modules.  In the situation at hand, we have
$$N:C_n(BG/Q)\lra C_n(G)^Q$$ of the form
$$[g_1|\cdots |g_n] \stackrel{N}{\mapsto} \sum_{q\in Q}[q(g_1)|\cdots |q(g_n)].$$
Since the domain is free abelian, $N$ is injective, but it is not surjective.  Denote the
quotient complex $\text{coker}(N)$ by $D_\bullet$.

\begin{prop}\label{quotientcomplex}  Suppose $Q=\zz/p$, where $p$ is a prime.  Then $D_0 = 0$, and for
$n>0$, $D_n = C_n(G^Q;\zz/p)$.  Consequently, there is an isomorphism
$$h_n(D_\bullet) \cong \tilde{H}_n(G^Q;\zz/p).$$
\end{prop}

\begin{proof} The map $N:C_0(BG/Q)\ra C_0(G)^Q$ is simply the identity map
$\zz\ra\zz$, whence the assertion for $D_0$.  If $n>0$, then since $p$ is prime,
the orbit of any $g\in G$ has length $1$ or $p$, and hence the same is true of
any $[g_1|\cdots |g_n]$.  The basis elements of $C_n(G)^Q$ are the various
$$\sum_{q\in Q}[q(g_1)|\cdots |q(g_n)]$$
together with the elements $[g_1|\cdots |g_n]$ for $g_i\in G^Q$.  If 
$[g_1|\cdots |g_n]$ is a generator of $C_n(BG/Q)$, then 
$$N:[g_1|\cdots |g_n]\mapsto \sum_{q\in Q}[q(g_1)|\cdots |q(g_n)].$$
If the orbit length of $[g_1|\cdots |g_n]$ is $p$, then we obtain the
corresponding basis element of $C_n(G)^Q$.  If, however, each $g_i\in G^Q$,
then
$$N:[g_1|\cdots |g_n]\mapsto p[g_1|\cdots |g_n].$$
Thus,
$$D_n = \frac{\zz\{[g_1|\cdots |g_n]: g_i\in G^Q\}}{p[g_1|\cdots |g_n]=0} \cong C_n(G^Q;\zz/p).$$
The boundary map in $D_\bullet$ is clearly that of $C_\bullet(G^Q;\zz/p)$; this proves
the final assertion.
\end{proof}

\begin{cor}\label{lesinvariant} If $Q=\zz/p$, then there is a long exact sequence
$$\ra H_n(BG/Q;\zz)\ra H_n^Q(G;\zz)\ra \tilde{H}_n(G^Q;\zz/p)\ra H_{n-1}(BG/Q;\zz)\ra.$$
\hfill $\qed$
\end{cor}

\begin{cor}\label{invariantquotient} Suppose that $Q=\zz/p$ and that $G^Q=\{e\}$.  Then there is an isomorphism
$$H_\bullet(BG/Q;\zz)\cong H_\bullet^Q(G;\zz).$$
\end{cor}

\begin{proof} In this case, we have $\tilde{H}_\bullet(G^Q;\zz/p) = 0$.
\end{proof}

For arbitrary $Q$, we also have the following.

\begin{prop} If $|Q|$ is invertible in $A$, then the map
$$H_\bullet(BG/Q;A)\lra H_\bullet^Q(G;A)$$ is an isomorphism.
\end{prop}

\begin{proof} The kernel and cokernel of the map
$$N:C_n(BG/Q;A)\lra C_n(G;A)^Q$$ are both annihilated by $|Q|$.  Since
$|Q|$ is invertible in $A$, we see that $N$ is an isomorphism.
\end{proof}

\section{Basic Calculations}\label{basic}
We begin our low-tech approach to calculation with a few bare-hands examples.

\begin{prop} $H_0^Q(G;A) = A$.
\end{prop}

\begin{proof} We have $C_0(G;A)^Q = A$ and the boundary map $$C_1(G;A)^Q\ra C_0(G;A)^Q$$ is the
zero map.
\end{proof}

The major difficulty in computing $H_\bullet^Q(G;A)$ is that the chain complex $C_\bullet(G;A)$
is so large (even if $G$ is finite).  Unfortunately, there is not much to be done about this.  Here
is an explicit example.

Let $\zz/2$ act on $\zz$ via $n\mapsto -n$.  Then $$C_i(\zz)^{\zz/2} = \zz\{[n_1|\cdots |n_i] + 
[-n_1|\cdots |-n_i],[0|\cdots |0]\}.$$

\begin{prop}\label{h1z} $H_1^{\zz/2}(\zz)\cong \zz/2$.
\end{prop}

\begin{proof}  We have
$$H_1^{\zz/2}(\zz) = \frac{\zz\{[n]+[-n],[0]: n>0\}}{\text{im}\{d:C_2(\zz)^{\zz/2}\ra C_1(\zz)^{\zz/2}\}}.$$
Define $f:C_1(\zz)^{\zz/2}\ra\zz/2$ by $f([n]+[-n]) = n \mod 2$ and $f([0]) = 0$.  Then $f$ is surjective.
Moreover, if $[n_1|n_2]+ [-n_1|-n_2] \in C_2(\zz)^{\zz/2}$, then
\begin{eqnarray*}
f(d([n_1|n_2] + [-n_1|-n_2])) & = & f(([n_1]+[-n_1]) + ([n_2]+[-n_2]) - \\
                              &   & {}\quad([n_1+n_2]+[-(n_1+n_2)])) \\
                              & = & n_1-(n_1+n_2) + n_2 \mod 2 \\
                              & = & 0 \mod 2
\end{eqnarray*}
and $f(d([0|0])) = f([0]-[0]+[0]) = 0$.  Thus, $f$ induces a surjection $\overline{f}:H_1^{\zz/2}(\zz)\ra\zz/2$.
It remains to show that $\overline{f}$ is injective.  For this it suffices to show
\begin{enumerate}
\item $[2k]+[-2k] = [0]$ in $H_1^{\zz/2}(\zz)$, and
\item $[2k+1]+[-(2k+1)] = [1]+[-1]$ in $H_1^{\zz/2}(\zz)$.
\end{enumerate}
Observe that $[0]$ is the identity element in $H_1^{\zz/2}(\zz)$ as $[0]=d([0|0])$.

Consider the relations
\begin{eqnarray*}
d([m|m]+[-m|-m]) & = & 2([m]+[-m]) - ([2m]+[-2m]) \\
d([m|-m]+[-m|m]) & = & 2([m]+[-m]) - 2[0].
\end{eqnarray*}
Note that the second equation implies that $H_1^{\zz/2}(\zz)$ is $2$-torsion.  Subtracting these shows
that $[2m]+[-2m] = 2[0]$ and hence $[2m]+[-2m] = 0$ in $H_1^{\zz/2}(\zz)$.

We also have the following equations
\begin{eqnarray*}
\lefteqn{d([m|m+3]+[-m|-(m+3)])  =}  \\
                  & &  ([m]+[-m]) + ([m+3]+[-(m+3)]) - ([2m+3]+[-(2m+3)]) \\
\lefteqn{d([m+1|m+2] + [-(m+1)|-(m+2)])  =} \\
                 & &  ([m+1]+[-(m+1)]) + ([m+2]+[-(m+2)]) \\
                 & & {}- ([2m+3]+[-(2m+3)])
\end{eqnarray*}
so that in $H_1^{\zz/2}(\zz)$ we have (by subtracting these)
\begin{eqnarray*}
\lefteqn{([m]+[-m]) + ([m+3]+[-(m+3)]) =  } \\
               & & {}\qquad\qquad([m+1]+[-(m+1)]) + ([m+2]+[-(m+2)]).
\end{eqnarray*}
Setting $m=0,2,4,\dots$ we obtain (since $[m]+[-m] = 0 = [m+2]+[-(m+2)]$)
$$[1]+[-1] = [3]+[-3] = \cdots = [2m+1]+[-(2m+1)] = \cdots.$$
This completes the proof.
\end{proof}

Now consider the resolution $F_\bullet\stackrel{\varepsilon}{\ra}\zz$ obtained
from $S^1\simeq B\zz$:
$$F_\bullet: \quad 0\lra\zz[t,t^{-1}]\stackrel{t-1}{\lra}\zz[t,t^{-1}]\stackrel{\varepsilon}{\lra}
\zz\lra 0.$$
To compute the complex $(F_\bullet\otimes_G\zz)^{\zz/2}$ of $\zz/2$-invariants, we need an
augmentation preserving chain map $\alpha_*:F_\bullet\ra F_\bullet$ commuting with the action.
The map $\alpha_0:t^n\mapsto t^{-n}$ and $\alpha_1:t^n\mapsto -t^{-(n+1)}$ does the job and
one finds that the complex $(F_\bullet\otimes_G\zz)^{\zz/2}$ is
$$0\lra 0\lra\zz.$$  The first homology of this complex vanishes, of course, and so we cannot
replace $C_\bullet(\zz)$ with $F_\bullet\otimes_G\zz$ to calculate $H_\bullet^{\zz/2}(\zz)$.

A moment's reflection reveals why this happens:  the invariants functor is left exact, but not
right exact, so there is no reason to believe that we can switch complexes at will.  From the
point of view of equivariant topology, there is a clearer answer:  the homotopy equivalence
$B\zz\ra S^1$ cannot be a $\zz/2$-equivariant equivalence.  Indeed, a map $f:X\ra Y$ of
$Q$-spaces is a $Q$-homotopy equivalence if the induced map $f^H:X^H\ra Y^H$ is a homotopy
equivalence for each $H\le Q$.  In the case of $B\zz\ra S^1$, $Q=\zz/2$, the $\zz/2$-fixed
points are $(B\zz)^{\zz/2}\simeq *$ while $(S^1)^{\zz/2} = S^0$.

Let us proceed with some general calculations.

\begin{prop}\label{invariantiso} Let $G$ be a group and let $A$ be an abelian group with trivial $G$-
and $Q$-actions.  Assume that $|Q|$ is
invertible in $A$.  Then the natural map
$$i_*:H_\bullet^Q(G;A)\lra H_\bullet(G;A)^Q$$ is an isomorphism.
\end{prop}

\begin{proof} This is a standard fact in equivariant topology, proved using the
transfer map.  We work directly with chains, however.  Let $[z]\in H_\ell^Q(G;A)$ and suppose
$i_*[z]=0$.  Then $z=d\tau$ for some $\tau\in C_{\ell+1}(G;A)$.  Since $z\in C_\ell(G;A)^Q$, we
have $qz=z$ for any $q\in Q$.  Thus, if the $Q$-orbit of $\tau$ has length $n$, we have
$$d\biggl(\sum_{q\in Q} q\tau\biggr) = \sum_{q\in Q} q(d\tau) = \sum_{q\in Q} qz = nz.$$
Write $\sigma = \sum_{q\in Q} q\tau$.  Then $q\sigma = \sigma$ for all $q\in Q$.  Since $n$ is
invertible in $A$ by assumption (note $n\mid |Q|$), we have $\frac{1}{n}\sigma
\in C_{\ell+1}(G;A)^Q$ and $d(\frac{1}{n}\sigma) = z$.  So $[z]=0$ and $i_*$ is injective.

To show that $i_*$ is surjective, let $[z]\in H_\ell(G;A)^Q$.  Then for each $q\in Q$, $qz=z+d\tau_q$
for some $\tau_q\in C_{\ell+1}(G;A)$.  Write $w=\sum_{q\in Q} qz$.  Then $w\in C_\ell(G;A)^Q$ and
$dw=0$.  We have
\begin{eqnarray*}
i_*[w] & = & \sum_{q\in Q} [qz] \\
       & = & \sum_{q\in Q} [z] \\
       & = & nz
\end{eqnarray*}
for some $n$ invertible in $A$.  So $i_*(\frac{1}{n}[w]) = [z]$ and $i_*$ is surjective.
\end{proof}

Note that as a corollary of the proof we have the following.

\begin{prop}\label{kerann} The subgroup $\text{\em ker}(i_*)$ is annihilated by $|Q|$.
\hfill $\qed$
\end{prop}

In particular, if $H_\ell(G;A)^Q=0$ for some $\ell$, we see that $H_\ell^Q(G;A)$ is annihilated
by the order of $Q$.  Thus, we have the following result.

\begin{prop}\label{hiz2tor} For each $i\ge 1$, the group $H_i^{\zz/2}(\zz)$ is $2$-torsion.
\hfill $\qed$
\end{prop}

Calculating $H_i^{\zz/2}(\zz)$ by hand from the bar resolution is essentially impossible, however.
We will complete the computation in Section \ref{calculate}.

\begin{prop}\label{divisible} If $G$ is a divisible abelian group,
then $$H_1^Q(G;\zz) = H_1(G;\zz)^Q=G^Q.$$
\end{prop}

\begin{proof}  We have the natural map
$$f_*:H_1(G^Q)\lra H_1^Q(G).$$ This map has a splitting $\tau:H_1^Q(G)\ra H_1(G^Q)$ defined by
$\tau([z_1]+\cdots [z_n]) = [z_1\cdots z_n]$, where $z_1,\dots ,z_n$ is a $Q$-orbit.  This is
well-defined since $$q(z_1\cdots z_n) = q(z_1)\cdots q(z_n) = z_1\cdots z_n$$ so that $z_1\cdots z_n\in G^Q$.

We now need the following relation in $H_1^Q(G)$:
\begin{equation}\label{h1rel}
[z_1^n]+\cdots + [z_n^n] = n[z_1\cdots z_n].
\end{equation}
Assuming this, we show that $f_*$ is surjective by noting that if $[w_1]+\cdots +[w_m]$ is an
orbit, write $w_1=z_1^m,\dots ,w_m=z_m^m$ and then
\begin{eqnarray*}
[w_1]+\cdots [w_m] & = & [z_1^m]+\cdots +[z_m^m] \\
                   & = & m[z_1\cdots z_m] \\
                   & = & f_*(m[z_1\cdots z_m]).
\end{eqnarray*}
But $f_*$ is injective since $\tau f_*([z]) = \tau[z] = [z]$.

It remains to show (\ref{h1rel}).  If $[z_1]+\cdots +[z_n]\in G_1(G)^Q$, set
$z_i=q_i(z_1)$, $i=2,\dots ,n$.  Consider the two relations
$$\sum_{j=1}^{m-1} ([z_1^j,z_i] + \cdots [z_m^j,z_m]) \stackrel{d}{\mapsto} m([z_1]+\cdots +[z_m]) - 
([z_1^m] +\cdots [z_m^m])$$
and
\begin{eqnarray*}
\lefteqn{\sum_{j=1}^{m-1} ([z_1\cdots z_j,z_{j+1}] + [q_2(z_1\cdots z_j),q_2(z_{j+1})] +\cdots +
[q_m(z_1\cdots z_j),q_m(z_{j+1})]) } \qquad\qquad\qquad\qquad\qquad\qquad\qquad\; \\
 & & {}\stackrel{d}{\mapsto} m([z_1]+\cdots +[z_m]) - m[z_1\cdots z_m].
 \end{eqnarray*}  
 Subtracting these yields (\ref{h1rel}).
 \end{proof}
 
 For example, if $G=\gm({\mathbb C})$ with
 $Q=\text{Gal}({\mathbb C}/{\mathbb R})$, then $H_1^Q(\gm({\mathbb C}))= \gm({\mathbb R})$.
 
 \section{Transfer}\label{transfer}
 Before we can make more calculations, we need to discuss transfer maps.  Recall that if $K$
 is a subgroup of finite index in $G$, then there is a transfer map
 $$\text{tr}:H_\bullet(G;A)\lra H_\bullet(K;A).$$ If $j_*:H_\bullet(K;A)\ra H_\bullet(G;A)$ is
 the map induced by the inclusion $j:K\ra G$, then the composite $j_*\circ \text{tr}$ is
 $(G:K)\text{id}$.  There is a double coset formula for the composite $\text{tr}\circ j_*$; we
 shall need only a special case of it in Section \ref{calculate}.  It will be recalled there.
 
 We want to know that we have transfer maps in our setting.  To this end, assume that $K$
 is a $Q$-stable subgroup of finite index in $G$ (that is, $qK=K$ for each $q\in Q$, but it
 is not necessarily the case that $K\subseteq G^Q$).
 For arbitrary $K$ and $G$, the transfer map is
 induced by a chain map $\tau:C_\bullet(G)\ra C_\bullet(K)$ defined by choosing a set $E$ of
 right coset representatives of $K$ and setting
 \begin{eqnarray*}
 \lefteqn{\tau[g_1|\cdots |g_n] =} & &  \\
   & &  \sum_{x_j\in E} [x_jg_1\overline{x_jg_1}^{-1}|\overline{x_jg_1}g_2\overline{x_jg_1g_2}^{-1}|
 \cdots | \overline{x_jg_1\cdots g_{n-1}}g_n\overline{x_jg_1\cdots g_n}^{-1}],
 \end{eqnarray*}
 where $\overline{a}$ is the representative of $Ka$ in $E$.
 
 In order for this to commute with the $Q$-action, we need a set of representatives $E$ such that
 $q(\overline{g}) = \overline{q(g)}$ for every $q\in Q$.  If we can do this, then since the homotopy
 from $j_*\circ\text{tr}$ to $(G:K)\text{id}$ maps orbits to orbits, the formula $j_*\circ\text{tr}=(G:K)\text{id}$
 will hold in this setting.  In the case $K=\{e\}$, we can do this
 since the coset representative of $g\in G$ is just $g$ itself and hence the transfer map exists
 on the chain level $C_\bullet(G;A)^Q\ra C_\bullet(\{e\};A)^Q$.  We therefore obtain the following.
 
 \begin{prop}\label{orderkill} If $G$ is finite, then $H_i^Q(G;A)$ is annihilated by $|G|$ for
 all $i>0$.
 \end{prop}
 
 \begin{proof} Let $j:\{e\}\ra G$ be the inclusion.  Then $j_*\circ \text{tr} = |G|\text{id}$, and
 $H_i^Q(\{e\};A) = 0$ for all $i>0$.
 \end{proof}
 
 Another example is the following.  Let $Q=\zz/2$ act on $\zz/2^sk$, $k$ odd, by $n\mapsto -n$.  Let
 $K$ be the subgroup $\langle k\rangle$.  The set $$E=\biggl\{0,1,2k-1,2,2k-2,\dots ,\frac{k-1}{2},
 \frac{k+1}{2}\biggr\}$$ is a set of coset representatives of $K$ in $\zz/2^sk$ and if $g$ denotes the nonidentity
 element of $\zz/2$, then for $1\le i\le \frac{k-1}{2}$, $g(\overline{i}) = g(i) = 2k-i$ and
 $\overline{g(i)} = \overline{2k-i} = 2k-i$.  Moreover, $g(\overline{2k-i}) = g(2k-i) = i$
 and $\overline{g(2k-i)} = \overline{i} = i$.  So we get a transfer map on the level of chains in
 this case as well:
 $$\text{tr}:H_\bullet^Q(\zz/2^sk;A)\lra H_\bullet^Q(K;A)$$
 and this satisfies $j_*\circ\text{tr}=k\text{id}$, where $j:K\ra\zz/2^sk$ is the inclusion.
 
 \section{Computations}\label{calculate}
 In this section we carry out some calculations.  We begin with the case $Q=\zz/2$ acting on $G=\zn$
 by $x\mapsto -x$.  In keeping with our desire to remain as low-tech as possible for as long as possible,
 we break this into cases.
 
 \subsection{$\zz/2$ acting on $\zn$, $n$ odd}  We claim the following:
 $$H_\ell^{\zz/2}(\zn) = \begin{cases}
                            \zz & \ell = 0 \\
                            \zn & \ell = 4k-1, k\ge 1 \\
                            0 & \text{otherwise}.
                         \end{cases}$$
 The case $\ell=0$ is clear.  If $\ell$ is even, then $H_\ell(\zn)=0$.  By Proposition \ref{kerann},
 we see that $H_\ell^{\zz/2}(\zn)$ is $2$-torsion.  On the other hand, $H_\ell^{\zz/2}(\zn)$ is 
 annihilated by $n$ (Proposition \ref{orderkill}) and since $n$ is odd, $H_\ell^{\zz/2}(\zn)$ must
 vanish for $\ell$ even.
 
 Now suppose that $\ell$ is odd and consider the map $i_*:H_\ell^{\zz/2}(\zn)\ra H_\ell(\zn)^{\zz/2}$.
 The kernel of this is $2$-torsion, but it is also $n$-torsion and hence $i_*$ is injective.
 
 If $\ell = 4k-3$, then the $\zz/2$-action on $H_\ell(\zn) = \zn$ is nontrivial: $x\mapsto -x$ (this
 is an easy exercise in group cohomology).  Since $n$ is odd, we have $H_{4k-3}(\zn)^{\zz/2} = 0$
 and since $i_*$ is injective, $H_{4k-3}^{\zz/2}(\zn) = 0$.
 
 It remains to treat the case $\ell = 4k-1$.  Here the $\zz/2$-action on $H_\ell(\zn)=\zn$ is trivial
 so that $H_\ell(\zn)^{\zz/2}= \zn$.  We must show that $i_*$ is surjective.  Let $z$ be a generating cycle of $H_\ell(\zn)$.
 Denote by $\overline{z}$ the $\zz/2$-translate of $z$.  Since $\zz/2$ acts trivially on $H_\ell(\zn)$, we
 have $\overline{z}=z+d\sigma$ for some $(\ell+1)$-chain $\sigma$.  Now, $z+\overline{z}$ is an element
 of $C_\ell(\zn)^{\zz/2}$ and $$i_*([z+\overline{z}])
 = [z]+[\overline{z}] = 2[z].$$ Since $n$ is odd, $2[z]$ also generates $H_\ell(\zn)$, and
 thus $i_*$ is an isomorphism.
 
 To summarize, then, we have proved that if $n$ is odd,
 $$i_*:H_\bullet^{\zz/2}(\zn;\zz) \stackrel{\cong}{\lra} H_\bullet(\zn;\zz)^{\zz/2}.$$

 \subsection{$\zz/2$ acting on $\zz/2k$, $k$ odd} Write $k=p_1^{i_1}\cdots p_m^{i_m}$. By Proposition \ref{invariantiso}, the
 map $$i_*:H_\bullet^{\zz/2}(\zz/2k;\zz/p^i) \lra H_\bullet(\zz/2k;\zz/p^i)^{\zz/2}$$
 is an isomorphism for each $p=p_j$, $i=i_j$.  The groups $H_\ell(\zz/2k;\zz/p^i)^{\zz/2}$
 are easily calculated:  they vanish if $p\ne p_j$ and
 $$H_\ell(\zz/2k;\zz/p_j^{i_j})^{\zz/2} = \begin{cases}
                                           \zz/p_j^{i_j} & \ell =4k-1, 4k \\
                                           0 & \text{otherwise}.
                                           \end{cases}$$
 
 It remains to calculate $H_\bullet^{\zz/2}(\zz/2k;\zz/2)$.  Let $K=\langle k\rangle$ be the subgroup
 of $\zz/2$-invariants.  Note that $(\zz/2k:K)=k$.  Let $j:K\ra \zz/2k$ be the inclusion.  We claim
 that
 $$j_*:H_\bullet^{\zz/2}(K;\zz/2)\lra H_\bullet^{\zz/2}(\zz/2k;\zz/2)$$
 is an isomorphism.  Indeed, since $\zz/2$ acts trivially on $K$, we have $$H_\bullet^{\zz/2}(K;\zz/2)
 = H_\bullet(K;\zz/2).$$  Consider the transfer map
 $$\text{tr}:H_\bullet^{\zz/2}(\zz/2k;\zz/2)\lra H_\bullet(K;\zz/2).$$  We have $j_*\circ \text{tr}
 = k\text{id} = \text{id}$.  Thus, $j_*$ is surjective.  To see that $j_*$ is injective, consider
 the composite
 $$H_\bullet(K;\zz/2)\stackrel{j_*}{\lra} H_\bullet^{\zz/2}(\zz/2k;\zz/2)\stackrel{i_*}{\lra}
 H_\bullet(\zz/2k;\zz/2)^{\zz/2}$$
 (note that the last is simply $H_\bullet(\zz/2k;\zz/2)$ since $\zz/2$ acts trivially on the
 mod 2 homology).  We claim that $i_*\circ j_*$ is an isomorphism.  Indeed, $i_*\circ j_*$ is the
 map induced by the inclusion $K\ra \zz/2k$.  Consider the commutative diagram
 $$\SMALL{\xymatrix{
 H_{2m+2}(\zz/2k;\zz/2)\ar[r] & H_{2m+1}(\zz/2k;\zz)\ar[r]^{\times 2} & H_{2m+1}(\zz/2k;\zz)\ar[r] & H_{2m+1}(\zz/2k;\zz/2) \\
 H_{2m+2}(K;\zz/2)\ar[u]\ar[r] & H_{2m+1}(K;\zz)\ar[u]\ar[r]^{\times 2} & H_{2m+1}(K;\zz)\ar[u]\ar[r] & H_{2m+1}(K;\zz/2)\ar[u]}}$$
 where the rows are the long exact sequences associated to the short exact coefficient sequence
 $$0\lra\zz\stackrel{\times 2}{\lra}\zz\lra\zz/2\lra 0.$$
 This diagram takes the form
 $$\xymatrix{
 \zz/2 \ar[r] & \zz/2k\ar[r]^{\times 2} & \zz/2k\ar[r]^{k\mapsto 1} & \zz/2 \\
 \zz/2\ar[u]^{i_*\circ j_*}\ar[r]^{\cong} & \zz/2\ar[u]^{1\mapsto k}\ar[r]^0 & \zz/2\ar[u]^{1\mapsto k}\ar[r]^\cong &
 \zz/2\ar[u]^{i_*\circ j_*}}$$
 An easy diagram chase shows that $i_*\circ j_*$ is an isomorphism.  Thus, $j_*$ is injective and hence an
 isomorphism.
 
 We can now piece this together to get a calculation of $H_\bullet^{\zz/2}(\zz/2k;\zz)$.  Consider the diagram
 of universal coefficient sequences
 $$\xymatrix{
 H_\ell(K)\otimes\zz/2\ar[d]\ar[r] & H_\ell(K;\zz/2)\ar[d]^\cong\ar[r] & \text{Tor}(H_{\ell-1}(K),\zz/2)\ar[d] \\
 H_{\ell}^{\zz/2}(\zz/2k)\otimes \zz/2\ar[r] & H_{\ell}^{\zz/2}(\zz/2k;\zz/2)\ar[r] & \text{Tor}(H_{\ell-1}^{\zz/2}(\zz/2k),\zz/2)}$$
 Note that $j_*:H_\ell(K)\ra H_\ell^{\zz/2}(\zz/2k)$ satisfies $$j_*\circ \text{tr} = k\text{id}_{H_\ell^{\zz/2}(\zz/2k)}.$$
 If $\ell\ge 2$ is even, then $H_\ell^{\zz/2}(\zz/2k;\zz)$ is 2-torsion since $H_\ell(\zz/2k;\zz)=0$.
 But then $k$ is invertible in $H_\ell^{\zz/2}(\zz/2k;\zz)$ so that $j_*\circ\text{tr}$ is an isomorphism.  But
 $H_\ell(K;\zz) = 0$ and hence $j_*\circ\text{tr}$ is the zero map.  This forces $H_{\ell}^{\zz/2}(\zz/2k;\zz)=0$.
 The above diagram then shows that the map
 $$H_\ell^{\zz/2}(\zz/2k;\zz/2)\lra\text{Tor}(H_{\ell-1}^{\zz/2}(\zz/2k),\zz/2)$$
 is an isomorphism and so if $r$ is odd, $H_r^{\zz/2}(\zz/2k;\zz)$ contains a single copy of $\zz/2$.  If
 $r=4s-3$, this is all $H_r^{\zz/2}(\zz/2k;\zz)$ contains (see the calculation above).
 If $r=4s-1$, then $H_r^{\zz/2}(\zz/2k;\zz)$ also contains a copy of $\zz/p^i$.  Thus
 $$H_{4s-3}^{\zz/2}(\zz/2k;\zz) = \zz/2$$
 and 
 $$H_{4s-1}^{\zz/2}(\zz/2k;\zz) = \zz/2\oplus \zz/p_i^{i_1} \oplus \cdots \oplus \zz/p_m^{i_m} = \zz/2k.$$
 Note that $\zz/2$ acts trivially on $H_{4s-1}(\zz/2k)$ and nontrivially on $H_{4s-3}(\zz/2k)$; the invariant
 subgroup of the latter is precisely $H_{4s-3}(K) = \zz/2$.  In summary, then, in the last two subsections
 we have proved the following.
 
 \begin{theorem}\label{nodd}  If $n\not\equiv 0\mod 4$, then
 $$i_*:H_\bullet^{\zz/2}(\zn;\zz)\lra H_\bullet(\zn;\zz)^{\zz/2}$$
 is an isomorphism. \hfill $\qed$
 \end{theorem} 
 
 \subsection{$\zz/2$ acting on $\zn$, $n\equiv 0\mod 4$}  The technique used in the previous section
 breaks down when $n\equiv 0 \mod 4$ because the subgroup of invariants has even index.  This then
 tells us that certain maps are zero, but not that they are isomorphisms.  Again, away from $p=2$,
 we have 
 $$i_*:H_\bullet^{\zz/2}(\zn;\zz/p^i)\stackrel{\cong}{\lra} H_\bullet(\zn;\zz/p^i)^{\zz/2}$$
 ($p^i\mid n$) and this is easily calculated.  Therefore we shall confine our attention to the
 calculation of $H_\bullet^{\zz/2}(\zn;\zz/2)$.
 
 Write $n=2^sm$, where $s\ge 2$ and $m$ is odd.  Denote by $M$ the subgroup $\langle m\rangle$.  Then
 $M\cong \zz/2^s$ and $(\zn:M)=m$.  Let $\alpha:M\ra\zn$ be the inclusion.  Then $\alpha_*\circ \text{tr}$
 is $m\text{id}$.  Taking $\zz/2$ coefficients, we see that
 $$\alpha_*:H_\bullet^{\zz/2}(M;\zz/2)\lra H_\bullet^{\zz/2}(\zn;\zz/2)$$
 is surjective.  We claim that $\alpha_*$ is injective as well.  For this we need a special case
 of the transfer double coset formula:  If $H$ is a normal subgroup of $G$, and $\varphi:H\ra G$ is
 the inclusion, then $$(\text{tr}\circ\alpha_*)(z) = \sum_{g\in G/H} gz,$$
 where $gz$ is obtained from $z$ via the conjugation action of $G$ on $H$.  In our case, everything in
 sight is abelian and so conjugation is trivial.  It follows that 
 $$(\text{tr}\circ\alpha_*)(z) = \sum_{g\in(\zn)/M} gz = mz.$$
 Since $m$ is odd, we have $\text{tr}\circ\alpha_* = \text{id}$.  Thus, $\alpha_*$ is injective
 as well.
 
 Thus, it remains to calculate $H_\bullet^{\zz/2}(\zz/2^s;\zz)$, and for this we need heavier
 machinery.  
 To this end, let $\Gamma=\zz/2$, $X=B(\zz/2^s)$, $Y=X^\Gamma=B(\langle 2^{s-1}\rangle)$.  We first
 calculate the homology of $X/\Gamma$.  If $p\ne 2$, we have $H_\bullet(X/\Gamma;\zz/p) = 
 H_\bullet(X;\zz/p)^{\zz/2}=0$.  We therefore devote our
 attention to the case of $\zz/2$ coefficients, and since homology and cohomology are dual over
 a field, we calculate cohomology instead. In what follows, unless otherwise specified, cohomology is
 with $\zz/2$ coefficients.
 
  Denote
 by $\pi$ the quotient map $X\ra X/\Gamma$.  Let $\beta:Y\ra X$ be the inclusion and note that
 $\pi\circ\beta$ maps $Y$ homeomorphically into $X/\Gamma$ as a closed subspace.  We first note
 the following fact:  $\beta^*:H^r(X)\ra H^r(Y)$ is the zero map for $r$ odd and the identity
 for $r$ even.  This is easily seen via a diagram chase using the long exact coefficient sequences
 {\SMALL $$\xymatrix{
 H^{2k-1}(X;\zz/2)\ar[r]\ar[d]^{\beta^*} & H^{2k}(X;\zz)\ar[r]^{\times 2}\ar[d]^{\beta^*} & H^{2k}(X;\zz)\ar[r]\ar[d]^{\beta^*} &
 H^{2k}(X;\zz/2)\ar[d]^{\beta^*} \\
 H^{2k-1}(Y;\zz/2)\ar[r]^{\cong} & H^{2k}(Y;\zz)\ar[r]^0 & H^{2k}(Y;\zz)\ar[r]^{\cong} & H^{2k}(Y;\zz/2) }$$}
 Denote by $\gamma$ the inclusion $\pi\circ\beta : Y\ra X/\Gamma$.  We need the following result.
 
 \begin{lemma} $\gamma^*:H^q(X/\Gamma)\lra H^q(Y)$ is the zero map for all $q>0$.
 \end{lemma}
 
 \begin{proof} Consider the commutative diagram obtained from the long exact sequences of the
 pairs $(X/\Gamma,Y)$ and $(X,Y)$:
 $$\xymatrix{
 H^{2k+1}(X/\Gamma)\ar[d]^{\pi^*}\ar[r]^{\gamma^*} & H^{2k+1}(Y)\ar[d]^{\cong}\ar[r] & H^{2k+2}(X/\Gamma,Y)\ar[d] \\
 H^{2k+1}(X)\ar[r]^0 & H^{2k+1}(Y)\ar[r]^{\cong} & H^{2k+2}(X,Y)}$$
 A diagram chase shows that $\gamma^*:H^{2k+1}(X/\Gamma)\ra H^{2k+1}(Y)$ is the zero map.
 
 To show $\gamma^*:H^{2k}(X/\Gamma)\ra H^{2k}(Y)$ is zero, we work directly with chains.  Since homology
 and cohomology are dual with $\zz/2$ coefficients, it suffices to show that $\gamma_*:H_{2k}(Y)\ra H_{2k}(X/\Gamma)$
 is zero.  By constructing a chain homotopy from the usual periodic resolution of $\zz$ over $\zz/2$ to
 the bar resolution, it is easy to see that
 $$z=\sum_{i_k=0}^1\sum_{i_{k-1}=0}^1\cdots \sum_{i_1=0}^1 [i_k\cdot 2^{s-1}|2^{s-1}|i_{k-1}\cdot 2^{s-1}|\cdots
 |i_1\cdot 2^{s-1}|2^{s-1}]$$
 is a generator for $H_{2k}(Y;\zz/2)$.  Let $\sigma_0=0|2^{s-1}$ and $\sigma_1=2^{s-1}|2^{s-1}$ and consider the
 $(2k+1)$-chain
 {\small $$w=\sum_{i_{k-1}=0}^1\cdots\sum_{i_1=0}^1 ([2^{s-2}|2^{s-2}|2^{s-1}|\sigma_{i_{k-1}}|\cdots |\sigma_{i_1}] + 
 [2^{s-2}|-2^{s-2}|2^{s-1}|\sigma_{i_{k-1}}|\cdots |\sigma_{i_1}])$$}
 in $C_\bullet(X/\Gamma;\zz/2)$.  A straightforward (if tedious) verification shows that $z=dw$ in $C_\bullet(X/\Gamma;\zz/2)$
 (one must keep in mind that, for example, $[2^{s-2}|2^{s-1}] = [-2^{s-2}|2^{s-1}]$ in $C_\bullet(X/\Gamma;\zz/2)$).
 Thus, $\gamma_*([z])=0$.
 \end{proof}
 
 Now, the lemma implies that the long exact sequence of the pair $(X/\Gamma,Y)$ breaks into
 short exact sequences
 \begin{equation}\label{xmodgamma}
 0\lra H^{q-1}(Y)\lra H^q(X/\Gamma,Y)\lra H^q(X/\Gamma)\lra 0.
 \end{equation}
 
 According to \cite{may}, p.~34, there is a long exact sequence (with $\zz/2$ coefficients)
 \begin{equation}\label{compseq}
 \tilde{H}^n((X/Y)/\Gamma)\ra H^n(X)\ra \tilde{H}^n((X/Y)/\Gamma)\oplus H^n(Y)\ra \tilde{H}^{n+1}((X/Y)/\Gamma).
 \end{equation}
 Note, however, that
$$(X/Y)/\Gamma = (X/\Gamma)/Y$$
and
$$\tilde{H}^\bullet((X/\Gamma)/Y) \cong H^\bullet(X/\Gamma,Y).$$
Thus, this sequence becomes
$$H^n(X/\Gamma,Y)\ra H^n(X)\ra H^n(X/\Gamma,Y)\oplus H^n(Y)\ra H^{n+1}(X/\Gamma,Y).$$
In this sequence, the map $H^n(X)\ra H^n(Y)$ is easily seen to be the usual restriction map.

 Since the map $H^{2k}(X)\ra H^{2k}(Y)$ is an isomorphism, the map
 $$H^{2k}(X)\lra H^{2k}(X/\Gamma,Y)\oplus H^{2k}(Y)$$ is injective.  We claim that the map
 $H^{2k+1}(X/\Gamma,Y)\ra H^{2k+1}(X)$ is surjective.  Indeed, consider the long exact sequences of
 the pairs $(X/\Gamma,Y)$ and $(X,Y)$:
 {\tiny$$\xymatrix{
 H^{2k}(Y)\ar[r]^(.4){0}\ar[d]^{\cong} & H^{2k+1}(X/\Gamma,Y)\ar[r]\ar[d]^{\pi^*} & H^{2k+1}(X/\Gamma)\ar[r]\ar[d]^{\pi^*} &
 H^{2k+1}(Y)\ar[r]^(.45){0}\ar[d]^{\cong} & H^{2k+2}(X/\Gamma,Y)\ar[d]^{\pi^*} \\
 H^{2k}(Y)\ar[r]^(.4){0} & H^{2k+1}(X,Y)\ar[r]^(.55){\cong} & H^{2k+1}(X)\ar[r]^{0} & H^{2k+1}(Y)\ar[r]^(.45){\cong} & H^{2k+2}(X,Y) }$$}
 The map $H^{2k+1}(X/\Gamma,Y)\ra H^{2k+1}(X)$ is the composite
 $$H^{2k+1}(X/\Gamma,Y) \stackrel{\pi^*}{\lra} H^{2k+1}(X,Y)\stackrel{\cong}{\lra} H^{2k+1}(X).$$
 The generator of $H^{2k+1}(X)$ is the function $f\in \text{Hom}(C_{2k+1}(X),\zz/2)$ given by
 $f(z)=1$, $f(w)=0$, $w\ne z$, where $z$ is the generating cycle
 $$z = \sum_{i_k=0}^{2^s-1}\sum_{i_{k-1}=0}^{2^s-1}\cdots\sum_{i_1=0}^{2^s-1} [1|{i_k}|1|\cdots |{i_1}|1]$$
 of $H_{2k+1}(X;\zz/2)$.  Define $h:C_{2k+1}(X/\Gamma)\ra\zz/2$ by
 $h(\{z\}) = 1$ and $h(\{w\})=0$ for $\{w\}\ne \{z\}$, where $\{z\}$ is the
 $\Gamma$-orbit of $z$.  Then $\pi^*(h) = f$ and so $\pi^*$ is surjective.
 
 From this, we conclude via (\ref{compseq}) that we have the following exact sequences for all $k\ge 0$
 {\small
 \begin{equation*}
 0\ra H^{2k}(X)\ra H^{2k}(X/\Gamma,Y)\oplus H^{2k}(Y)\ra H^{2k+1}(X/\Gamma,Y)\ra H^{2k+1}(X)\ra 0;
 \end{equation*}
 }
 and
 $$0\ra H^{2k+1}(X/\Gamma,Y)\oplus H^{2k+1}(Y)\ra H^{2k+1}(X/\Gamma,Y)\ra 0.$$
 Thus, $H^0(X/\Gamma,Y) = 0$, $H^1(X/\Gamma,Y) \cong \zz/2$, $H^2(X/\Gamma,Y) \cong (\zz/2)^2$,$\dots$,
 $H^q(X/\Gamma,Y) \cong (\zz/2)^q$.  Substituting this into (\ref{xmodgamma}) yields the following
 result.
 
 \begin{theorem}\label{n0mod4xmodgamma}  If $s\ge 2$, then $H_1(X/\Gamma;\zz/2)\cong \zz/2$ and for all $i\ge 2$
 $$H_i(X/\Gamma;\zz/2)\cong (\zz/2)^{i-1}.$$
 \hfill $\qed$
 \end{theorem}

Now, we need to compute the integral homology $H_\bullet(X/\Gamma;\zz)$.  We first note the following fact.

\begin{prop} If $\ell\ne 4k-1$, then $H_\ell(X/\Gamma;\zz)$ is $2$-torsion.
\end{prop}

\begin{proof} Note that we have a transfer map $$\mu:H_\bullet(X/\Gamma;\zz)\lra H_\bullet(X;\zz)$$
satisfying the following
$$\pi_*\circ\mu(\alpha) = 2\alpha\qquad \mu\circ\pi_*(\beta) = \beta + \overline{\beta}$$
where $\overline{\beta}$ is the image of $\beta$ under the action of $\zz/2$ (see \cite{borel}, p.~37).
If $\ell$ is even, then since $H_\ell(X;\zz)=0$, we see that $2\cdot H_{\ell}(X/\Gamma;\zz)=0$.

If $\ell=4k-3$, then $\mu([z]) = [z]+[\overline{z}] = 0$ (since $\overline{\beta}=-\beta$ in $H_\ell(X;\zz)$), 
but $\pi_*(\mu([z]))=2[z]$ so that $2[z]=0$, unless $z$ is actually a chain over $Y=X^{\Gamma}$.  In that
case, $\mu([z])=[z]$, but $2[z]=0$ anyway since $[z]$ lies in the image of $H_\ell(Y;\zz)=\zz/2$.
\end{proof}

\begin{prop} $H_{4k-1}(X/\Gamma;\zz)$ contains a single copy of $\zz/2^s$.
\end{prop}

\begin{proof} Let $[z]$ be any generator of $H_{4k-1}(X;\zz)\cong\zz/2^s$.  Let $[w]$ be a
generator of $H_{4k-1}(Y;\zz)\cong \zz/2$.  Then $2^{s-1}[z]=[w]$ and
\begin{eqnarray*}
\mu\circ\pi_*(2^{s-1}[z]) & = & \mu\circ\pi_*([w]) \\
                                & = & \mu([w]) \\
                                & = & [w] \\
                                & \ne & 0.
\end{eqnarray*}
Thus, $\pi_*(2^{s-1}[z])\ne 0$ and hence $2^{s-1}[z]\ne 0$ in $H_{4k-1}(X/\Gamma;\zz)$.  On the other
hand,
\begin{eqnarray*}
2^s[z] & = & \pi_*\circ\mu(2^{s-1}[z]) \\
       & = & \pi_*(2^{s-1}[z]+2^{s-1}[\overline{z}]) \\
       & = & \pi_*(2^s[z]) (\text{since}\,[\overline{z}]=[z]) \\
       & = & \pi_*(0) \\
       & = & 0.
\end{eqnarray*}
Thus, every generator of $H_{4k-1}(X;\zz)$ maps to an element of order $2^s$ in $H_{4k-1}(X/\Gamma;\zz)$
and since these generators are all multiples of each other in $H_{4k-1}(X;\zz)$, the same is true in
$H_{4k-1}(X/\Gamma;\zz)$.
\end{proof}

Using the previous two propositions and the universal coefficient theorem, we then obtain the following
result.

\begin{theorem}\label{xmodgammaz}  For all $\ell>0$,
$$H_\ell(X/\Gamma;\zz) \cong \begin{cases}
                              (\zz/2)^{k-1} & \ell=2k \\
                              (\zz/2)^{2k-1} & \ell=4k-3 \\
                              \zz/2^s \oplus (\zz/2)^{2k-1} & \ell=4k-1.
                              \end{cases}$$
\hfill $\qed$
\end{theorem}

Now to complete the calculation of $H_\bullet^{\zz/2}(\zz/2^s;\zz)$, we use Corollary \ref{lesinvariant}.  There
is a long exact sequence
$$\ra H_n(X/\Gamma;\zz)\ra H_n^{\zz/2}(\zz/2^s;\zz)\stackrel{p_*}{\ra} \tilde{H}_n(\zz/2;\zz/2) \stackrel{\partial}{\ra} H_{n-1}(X/\Gamma;\zz)\ra.$$
We claim that the map $p_*$ is surjective.  Indeed, let $z$ be a generating cycle for $H_n(\zz/2;\zz/2)$.  If $n$ is
odd, then $z$ is an integral cycle and $z\in Z_n(\zz/2^s;\zz)^{\zz/2}$ so that $p_*[z]=[z]$.  If $n$ is even, then
$dz=2w$ where $w$ generates $H_{n-1}(\zz/2;\zz)$.  But $2[w]=0$ in $H_{n-1}(X/\Gamma;\zz)$ (since $2[w]=0$ in $H_{n-1}(X;\zz)$) and so $\partial[z]=[2w]=0$.
It follows that $p_*$ is surjective in this case as well.

Thus, for each $q\ge 1$, we have a short exact sequence
$$0\ra H_q(X/\Gamma;\zz)\ra H_q^{\zz/2}(\zz/2^s;\zz)\ra H_q(\zz/2;\zz/2)\ra 0.$$
This clearly splits for $q\ne 4k-1$ since each term is $2$-torsion.  When $q=4k-1$ this still splits
since $H_{4k-1}(X/\Gamma;\zz)$ contains a copy of $\zz/2^s$ together with some copies of $\zz/2$, while
Proposition \ref{orderkill} implies that $H_{4k-1}^{\zz/2}(\zz/2^s;\zz)$ is annihilated by $2^s$ and hence there are no copies
of $\zz/2^{s+1}$ to provide a nonsplit quotient.  We therefore obtain the following calculation.

\begin{theorem}\label{n0mod4} For all $\ell>0$,
$$H_\ell^{\zz/2}(\zz/2^s;\zz) \cong \begin{cases}
                              (\zz/2)^{k} & \ell=2k \\
                              (\zz/2)^{2k} & \ell=4k-3 \\
                              \zz/2^s \oplus (\zz/2)^{2k} & \ell=4k-1.
                              \end{cases}$$
\hfill $\qed$
\end{theorem}

The interested reader may now use this to calculate $H_\bullet^{\zz/2}(\zn;\zz)$ for any $n\equiv 0\mod 4$.

\subsection{$\zz/2$ acting on $\zz$ via $n\mapsto -n$}  We are now able to calculate $H_\bullet^{\zz/2}(\zz)$.
We know this is all 2-torsion in positive degrees (Proposition \ref{hiz2tor}).  Set $\Gamma=\zz/2$, $X=B\zz$,
$Y=X^\Gamma\simeq *$.  The long exact sequence of the pair $(X/\Gamma, Y)$ (with $\zz/2$ coefficients)
$$H^q(X/\Gamma,Y)\lra H^q(X/\Gamma)\lra H^q(Y)\lra H^{q+1}(X/\Gamma,Y)$$
shows that $H^q(X/\Gamma,Y)\stackrel{\cong}{\lra} H^q(X/\Gamma)$ for all $q>0$.  Again, using \cite{may}, p.~34,
 we obtain a long exact sequence
$$H^q(X/\Gamma,Y)\ra H^q(X)\ra H^q(X/\Gamma,Y)\oplus H^q(Y)\ra H^{q+1}(X/\Gamma,Y).$$
Since $H^q(X)=0$ for $q\ge 2$ and $H^q(Y)=0$ for $q\ge 1$, we see that for $q\ge 2$
$$H^q(X/\Gamma,Y)\stackrel{\cong}{\lra} H^{q+1}(X/\Gamma,Y).$$
Also, $H^0(X/\Gamma,Y) = 0$ and we have an exact sequence
$$0\ra H^1(X/\Gamma,Y)\ra H^1(X)\ra H^1(X/\Gamma,Y)\ra H^2(X/\Gamma,Y)\ra 0.$$
If $H^1(X/\Gamma,Y) = 0$, then we would have the exact sequence
$$0\lra\zz/2\lra 0,$$
which is absurd.  So $H^1(X/\Gamma,Y) \cong\zz/2$ and the map $$H^1(X/\Gamma,Y) \ra H^2(X/\Gamma,Y)$$ is an isomorphism.
Thus, $H^q(X/\Gamma,Y)\cong\zz/2$ for all $q\ge 1$ and so $$H^q(X/\Gamma;\zz/2)\cong\zz/2$$ for all $q\ge 1$. 

\begin{theorem}\label{hiz}
$$\tilde{H}_i^{\zz/2}(\zz;\zz) \cong \begin{cases}
                                      \zz/2 & i\;\text{\em odd}\\
                                        0 & i\;\text{\em even}.
                                        \end{cases}$$
\end{theorem}

\begin{proof}
We know that $H_i^{\zz/2}(\zz;\zz)$ is 2-torsion for $i\ge 1$.  According to Corollary \ref{invariantquotient}, there is
an isomorphism $$H_i(X/\Gamma;\zz)\cong H_i^{\zz/2}(\zz;\zz).$$
The universal coefficient theorem applied to the homology of $X/\Gamma$ then yields the desired result.
\end{proof}


\begin{thebibliography}{99}

\bibitem{borel} A.~Borel, et.~al., {\em Seminar on transformation groups}, Annals
of Mathematics Studies {\bf 46}, Princeton University Press, 1960.

\bibitem{brown} K.~S.~Brown, {\em Cohomology of groups}, Springer--Verlag Graduate Texts
in Mathematics {\bf 87}, 1982.

\bibitem{knudsonwalker} K.~Knudson, M.~Walker, {\em Homology of linear groups via cycles in $BG\times X$}, preprint, 2003.

\bibitem{may} J.~P. May, {\em Equivariant homotopy and cohomology theory} (with contributions by
M.~Cole, G.~Comenza\~na, S.~Costenoble, A.~D.~Elmendorf, J.~P.~C.~Greenlees, L.~G.~Lewis, R.~J.~Piacenza,
G.~Triantafillou, and S.~Waner), CBMS Regional Conference Series in Mathematics {\bf 91}, American
Mathematical Society, Providence, 1996.



\end{thebibliography}
\end{document}